\documentstyle[amssymb,amsfonts]{amsart}

\def\Var{\hbox{\bf Var}}

\def\R{{\hbox{\bf R}}}

\def\P{{\hbox{\bf P}}}
\def\E{{\hbox{\bf E}}}

\font \roman = cmr10 at 10 true pt

\def\bv{{\hbox{\bf v}}}

\def\beq{\begin{equation}}
\def\eeq{\end{equation}}
\def\barray{\begin{eqnarray*}}
\def\earray{\end{eqnarray*}}
\def\be#1{ \begin{equation}\label{#1} }

\def\bas{\begin{align*}}
\def\eas{\end{align*}}
\def\bi{\begin{itemize}}
\def\ei{\end{itemize}}

\def\dist{{\hbox{\roman dist}}}

\def\Z{{\hbox{\bf Z}}}

\newenvironment{proof}{\noindent {\bf Proof} }{\endprf\par}
\def \endprf{\hfill  {\vrule height6pt width6pt depth0pt}\medskip}
\def\emph#1{{\it #1}}
\def\textbf#1{{\bf #1}}

\def\Bv{{\mathbf v}}

\def\BE{{\mathbf E}}

\def\BP{{\mathbf P}}

\def\BZ{{\mathbf Z}}

\def\CP{{\mathcal P}}

\def\ep{{\epsilon}}

\theoremstyle{plain}
  \newtheorem{theorem}[subsection]{Theorem}
  \newtheorem{conjecture}[subsection]{Conjecture}

  \newtheorem{fact}[subsection]{Fact}
  \newtheorem{lemma}[subsection]{Lemma}

\theoremstyle{remark}
  \newtheorem{remark}[subsection]{Remark}

\theoremstyle{definition}

\include{psfig}

\begin{document}

\title{ Random Discrete Matrices }

\author{  Van Vu}
\address{Department of Mathematics, Rutgers, Piscataway, NJ 08854}
\email{vanvu@@math.rutgers.edu}

\thanks{V. Vu is an A. Sloan  Fellow and is supported by  NSF Career Grant 0635606.}

\begin{abstract}
In this survey, we discuss some basic problems concerning random
matrices with discrete distributions. Several new results, tools
and conjectures will be presented.
\end{abstract}

\maketitle
\section{Introduction}
Random matrices is an important area of mathematics, with strong
 connections to many other areas (mathematical physics, combinatorics,
 theoretical computer science, to mention a few).

There are two types of random matrices: continuous and discrete. The
continuous models have an established theory (see \cite{Meh}, for
instance). On the other hand, the discrete models are still not very
well understood. In this survey, we discuss a few basic problems
concerning these models. The topics to be discussed are:
\begin{itemize}

\item The limiting  distribution of the spectrum (Section
\ref{limit}).

\item The spectral norm and the second largest eigenvalue
(Sections \ref{norm}, \ref{norm2}).

\item  Determinant (Section \ref{det}).

\item Rank and Singular probability (Sections \ref{rank},
\ref{rank2}).

\item The condition number (Section \ref{cond}).

\item Tools from additive combinatorics (Sections \ref{LO},
\ref{randomwalk}, \ref{inverseLO}).

\end{itemize}

 {\bf Notations.} We denote by $M_n$ the $n$ by $n$ random matrix
whose entries are i.i.d Bernoulli random variables (taking values
$1$ and $-1$ with probability $1/2$). This matrix is not symmetric.
Symmetric matrices often come from graphs.  We denote by $Q(n,p)$
the adjacency matrix of the Erd\H {o}s-R\'enyi random graph
$G(n,p)$. Thus $Q(n,p)$ is a random symmetric matrix whose upper
diagonal entries are i.i.d random variables taking value 1 with
probability $p$ and $0$ with probability $q=1-p$. Another popular
model for random graphs is that of random regular graphs. A random
regular graph $G_{n,d}$ is obtained by sampling uniformly over the
set of all simple $d$-regular graphs on the vertex set $\{1, \dots,
n \}$. The adjacency matrix of this graph is denoted by $Q_{n,d}$.

In the whole paper, we assume that $n$ is large. The asymptotic
notation is used under the assumption that $n \rightarrow \infty$.
We write $A \ll B$ if $A= o(B)$.  $c$ denotes a universal
constant. All logarithms have natural base, if not specified
otherwise.

\section {The universality principle}

Intuitively, one would expect a universal behavior among random
models of the same object. For random matrices in particular,  one
would expect the distributions of specific eigenvalues be the same
(after a proper normalization), regardless the model. Thus, given a
theorem for continuous models, it is often simple to come up with a
reasonable conjecture for discrete ones. For instance, there are
fairly  accurate tail estimates for the smallest singular value of a
random matrix whose entries are i.i.d Gaussians (see for example
Theorem \ref{theo:condition11} in Section \ref{cond}). It would be
natural to try to prove similar estimates for a random matrix whose
entries are i.i.d Bernoulli. However, this kind of task is usually
challenging, as the tools developed for  continuous models are
typically  not applicable in a discrete setting. In the  last few
sections (Sections \ref{LO}, \ref{randomwalk}, \ref{inverseLO}) of
this survey we will present new tools developed recently in order to
treat the discrete models. These tools, among others, reveal an
intriguing connection between the theory of random matrices and
additive combinatorics.

For random graphs, there is a specific conjecture  which establishes
the universality between the two models $G(n,p)$ and $G_{n,d}$
(Erd\H {o}s-R\'enyi graphs and random regular graphs).

\begin{conjecture} \label{conjecture1} (Sandwich Conjecture) \cite{KimV} For $ d \gg
\log n$, there is a joint distribution (or coupling) on random
graphs $H, G_{n,d}, G$ such that

\begin{itemize}

\item    $H$ has the same distribution as $G(n,p_1)$ where  $p_1=
\frac{d}{n}(1 - \frac{c\sqrt {d \log n}}{n} )$
  and $G$ has the same distribution as $G(n,p_2)$  where
  $p_2= \frac{d}{n}(1+ \frac{c\sqrt {d \log n}}{n}) $.

\item  $ \P( H\subset G_d ) = 1- o(1). $

\item  $ \P( G_{n,d} \subset G) = 1- o(1). $
\end{itemize}

\end{conjecture}

The conjecture asserts that a random regular graph can be
approximated from both below and above by  Erd\H {o}s-R\'enyi graphs
of approximately the same densities. The conjecture has been proved
for $d \le n^{1/3-o(1)}$ \cite{KimV}.

\begin{theorem} \label{theorem:sandwich}  The sandwich conjecture holds for $ \log n \ll d
\ll n^{1/3}/\log^2 n $.

\end{theorem}

The main difficulty when dealing with the random regular graph
$G_{n,d}$ is that the (upper diagonal) entries of its adjacency
matrix are not independent variables. But using  Theorem
\ref{theorem:sandwich}, one can often deduce information  about
the spectrum of $G_{n,d}$ using information about the spectrum of
$G(n,p)$.

\section {Limiting distribution} \label{limit}

One of the cornerstones of the theory of random matrices is Wigner's
semi-circle law, which established the limiting distribution of a
certain class of random symmetric matrices \cite{Wig}. We present
here a more general version, due to Arnold \cite{Ar}.

Let $a_{ij}$, $1 \le i \le j \le n$, be i.i.d random variables with
common variance one and distribution function $F(x)$ such that
$\int_{0}^{\infty} |x|^k dF < \infty $ for all $k=1,2, \dots$.  Let
$A_n$ be the random symmetric matrix of size $n$ whose upper
diagonal entries are $\xi_{ij} = a_{ij}/ 2 \sqrt n$. Let $\lambda_1
\ge \dots \ge \lambda_n$ be the (real) eigenvalues of $A_n$. Define

$$W_n(x) := \frac{1}{n} |\{i| \lambda_i \le  x \}| . $$

Let $W(x)$ denote the semi-circle density function

$$W(x) := \frac{2}{\pi} \sqrt {1 -x^2} $$

\noindent for $|x| \le 1$ and $W(x):=0$ otherwise.

\begin{theorem} \label{theorem:Wigner} (Semi-circle law) With
probability one,
$$\lim_{n \rightarrow \infty} W_n(x) = W(x) . $$
\end{theorem}

In order to prove the semi-circle law, Wigner introduced the
so-called trace method, the heart of which is the calculation of
the expectation of $Trace (A_n^k)$ for $k=1,2, \dots.$ This method
is useful for many other problems (see Section \ref{norm} for
example).

Let us now turn to the special matrix $Q(n,p)$. The entries of
$Q(n,p)$ have variance $\sigma^2= p(1-p)$.   Dividing each entry of
$Q(n,p)$ by $ \sigma $, we obtain a matrix $Q'(n,p)$ whose entries
have common variance one. However, one cannot apply Theorem
\ref{theorem:Wigner} directly as the entries of $Q'(n,p)$ do not
have bounded moments when $p$ tends to zero with $n$. On the other
hand, by applying  Wigner's trace method, one can prove the
following theorem

\begin{theorem} \cite{FK,VW} \label{theorem:randomgraph} There is a constant $c$ such that the following holds.
Let $\lambda_1 \ge \dots
\ge \lambda_n$ be the eigenvalues of $Q(n,p)$ where $p \ge n^{-1}
\log^c n $ and define

$$R^1_n(x):= \frac{1}{n} |\{i| \lambda_i \le  2 x\sqrt {np(1-p)}  \}|. $$

Then  with probability one,
$$\lim_{n \rightarrow \infty} R_n^1(x) = W(x) . $$
\end{theorem}

A corollary of a general theorem by Guionet and Zeitouni \cite{GZ}
shows that $R^1_n (x)$ and many other quantities concerning the
spectrum of $Q(n,p)$ are highly concentrated.

Next we discuss the situation with  the random regular graph
$G_{n,d}$. Define

$$R^2_n(x):=  \frac{1}{n} |\{i| \lambda_i \le  2 x\sqrt {d-1}   \}|  $$

\noindent and

$$W(d,x) : = \frac{d^2-d}{d^2 -4(d-1)x^2} \frac{2}{\pi} \sqrt {1-x^2}  $$

\noindent for $|x| \le 1$ and $W(d,x):=0$ otherwise. A theorem of
McKay \cite{Mc} on the spectrum of regular graphs (not necessarily
random) implies

\begin{theorem}  \label{theorem:McKay}
(Distribution of the eigenvalues in random regular graphs with fixed
degree) For  any fixed $d$ the following holds with probability one

$$\lim_{n \rightarrow \infty} R^2_n(x) = W(d,x) . $$
\end{theorem}

Observe that the limiting distribution $W(d,x)$ in this theorem is
not semi-circular because of the extra term $ \frac{d^2-d}{d^2
-4(d-1)x^2}$. On the other hand, it becomes arbitrarily close to the
semi-circle distribution if  $d$ is sufficiently large. Thus it is
reasonable to conjecture that if $d$ tends to infinity with $n$,
$Q_{n,d}$ follows the semi-circle law. However, McKay's proof used
Wigner's trace method and relied on the crucial fact that the graph
has few small cycles. Theorem \ref{theorem:McKay} still holds for
$d= n^{o(1)}$. But for  $d= n^{\ep}$ with any constant $\ep
>0$ the graph has too many small cycles and it seems very hard to apply this
method. On the other hand, using the sandwiching theorem (Theorem
\ref{theorem:sandwich}), Vu and Wu \cite{VW} proved that if $\log n
\ll d \ll n^{1/3}/\log ^2 n$ then with probability one

 $$\lim_{n \rightarrow \infty} R^2_n(x) = W(x) . $$

If the sandwich conjecture holds for all $d \gg \log n$, then this
statement can be extended for all $d \gg \log n$. Recently, Zeitouni
(private communication)
 suggested to the author another approach that also seems to work
for a wide range of $d$. The details will appear in \cite{VW}. For
results concerning more general models of random graphs, see
\cite{CLV, VW}.

To conclude this section, let us briefly discuss the case when $A_n$
is not symmetric. In this case, the eigenvalues are complex numbers.
Let $a_{ij}$, $1 \le i, j \le n$, be i.i.d complex random variables
with mean zero and variance one. Let $A_n$ be the random matrix with
entries $\xi_{ij} = a_{ij}/ \sqrt n$ and let $\lambda_1, \dots,
\lambda_n$ be the eigenvalues of $A_n$. Consider the two-dimensional
empirical distribution

$$\mu_n (x,y) := \frac{1}{n} |\{i| Re (\lambda_i) \le x, Im
(\lambda_i) \le y \}|. $$

It was proved by Girko \cite{Girko} and Bai \cite{Bai} that under
some weak conditions, $\mu_n (x,y)$ tends to the uniform
distribution over the unit disc.

\begin{theorem}  \label{theorem:GB} (Circular law) Assume that the $\xi_{ij}$ have finite sixth moment
and the joint distribution of the real and  imaginary part has
bounded density. Then with probability one $\mu_n (x,y)$ tends to
the uniform distribution over the unit disc in $\R^2$.
\end{theorem}

\section {The spectral norm} \label{norm}

The spectral norm of an $n$ by $n$ matrix  $A$ is defined as

$$\|A \| =\sup_{v \in \R^n, \|v\| =1} |Av|. $$

\noindent If $A$ is symmetric, then $\|A\|$ is the largest
eigenvalue of $A$ (in absolute value).

We consider the following general model of random symmetric
matrices. Let $\xi_{ij}$, $1\le i \le j \le n$, be independent
(but not necessarily identical) random variables with the
following properties

\begin{itemize}

\item $|\xi_{ij}| \le K$ for all $1 \le i \le j \le n$.

\item $\BE (\xi_{ij})=0$, for all $1 \le i < j \le n$.

\item $\Var (\xi_{ij}) = \sigma^2$, for all $1 \le i < j \le n$.

\end{itemize}

\noindent  For a moment, let us assume  that $\sigma$ and $K$ are
positive constants.

\vskip2mm

Define $\xi_{ji} = \xi_{ij}$ and consider the symmetric random
matrix $A_n=(\xi_{ij})_1^n$. Notice that for any matrix $A$,

$$||A||=\lim_{k \rightarrow \infty} Tr(A^k)^{1/k} . $$

This suggests that the trace method would be an effective tool for
bounding $\|A_n\|$. Indeed, all upper bounds mentioned below are
based on this method.  F\"uredi and Koml\'os \cite{FK} proved that
a.s.

$$
\| A_n \| \le 2 \sigma \sqrt n + cn^{1/3} \ln n. $$

The error term $cn^{1/3} \ln n$  was recently improved to $c
n^{1/4} \ln n$ by Vu \cite{Vunorm}. From below, Alon, Krivelevich
and Vu \cite{AKV} showed that a.s.

$$2 \sigma \sqrt n - c \ln n \le
\| A_n \|, $$

\noindent for some constant $c$. Putting these bounds together, we
have

\begin{theorem} \label{theorem:Vspectralnorm}  For a
 random matrix $A$ as above there is a positive constant
 $c=c(\delta, K)$  such that

$$  2 \sigma \sqrt n - c \ln n \le \| A \|  \le  2 \sigma \sqrt n + cn^{1/4} \ln n, $$

\noindent holds almost surely.

\end{theorem}

 In many situations
$\sigma$ and $K$ may depend on $n$. A typical example is when $A$ is
the ''normalized" adjacency matrix of $G(n,p)$ ($1$  and $0$ are
replaced by  $1-p$ and $-p$, respectively; this forces all entries
to have mean  zero) where $p$ is decreasing with $n$ ($p= n^{-\ep}$,
say). In this case, by following the proof of the upper bound in the
previous theorem, one can obtain

\begin{theorem} \label{theorem:Vspectralgeneral}
\cite{Vunorm} There are constants $c$ and $c'$ such that the
following holds.
 Let $\xi_{ij}$, $1\le i \le j \le n$ be independent
random variables, each of which has mean 0 and variance at most
$\sigma^2$ and is bounded in absolute value by $K$, where $\sigma
\ge c' n^{-1/2}K \ln^2 n$. Then almost surely

$$\|A_n \| \le 2\sigma \sqrt n + c (K\sigma)^{1/2} n^{1/4} \ln n. $$

\end{theorem}

One can also obtain a somewhat weaker bound by following the  proof
from \cite{FK}.

If one assumes more about the distributions of the entries
$\xi_{ij}$, one can obtain a sharper bound. Soshnikov and Sinai
\cite{SS} proved the following

\begin{theorem} \label{theorem:Sos}  (Spectral bound for random matrices with symmetric entries)
Let $A_n$ be a random symmetric matrix whose upper diagonal
entries $\xi_{ij}, 1\le i \le j \le n$ are independent random
variables satisfying

\begin{itemize}
\item $\xi_{ij}$ have symmetric distribution.

\item $\E( \xi_{ij}^2) = 1$ and $\E(\xi_{ii}^2) =O(1)$.

\item For all $m \ge 1$, $\E (\xi_{ij}^{2m}) = O(m)^m $.

\end{itemize}

Then a.s. $\| A_n \| = 2 \sqrt {n} + O( n^{-1/6})$.

\end{theorem}

In certain algorithmic applications, it is useful to have a tail
distribution for the spectral norm (see \cite{Optas, KV}, for
instance). Using Talagrand's inequality, one can show

\begin{theorem} \label{theorem:KV} \cite{KV, AKV} (Concentration of the spectral
norm) There is a positive constant
 $c=c(K)$ such that for any $t >0$

 $$\BP (\Big| \|A_n \|  -\E (\|A_n \|) \Big |  \ge  ct) \le 4e^{-t^2/32} . $$
 \end{theorem}

Similar results hold for  larger classes of random matrices  and
also for  other eigenvalues  (see \cite{AKV, Mec}).

     In the case of $Q(n,p)$, if $p$ is sufficiently large, then all
rows have a.s. roughly $np$ ones and the norm of $Q(n,p)$ is
(a.s.) $(1+o(1)) np$. But if $p$ is relatively small, this is no
longer  true. Krivelevich and Sudakov \cite{KSud} proved that
$\|Q(n,p)\|$ is almost surely

$$(1+o(1)) \max \{np, \sqrt D \} $$

\noindent where $D$ denotes the maximum degree of the (random)
graph. See \cite{SudS, Janson} for more results of this type.

\section{The second eigenvalue of random regular graphs}
\label{norm2}

Let $G$ be a graph on $n$ points and $A$  its adjacency matrix.
Let $\lambda_1 \ge \lambda_2 \ge \dots \ge \lambda_n$ be the
eigenvalues of $A$. If $G$ is $d$-regular, then $\lambda_1=d$. In
this case, a critical parameter of the graph is

$$\lambda (G):= \max \{|\lambda_2|, |\lambda_n| \}. $$

In the literature, $\lambda (G)$ is frequently called the second
eigenvalue of $G$. (This name is inaccurate but somewhat
convenient.) The good  way to think of $\lambda (G)$ is

$$\lambda(G) = \max_{\|v\|=1, v \cdot {\bf 1} =0} |v^{T} \cdot A v|,
$$

\noindent where ${\bf 1}$ is the all ones vector. One can also think
of $\lambda(G)$ as the spectral norm of the ''normalized" adjacency
matrix of $G$, where $1$ and $0$ are replaced by $(n-d)/n$ and
$-d/n$, respectively.

One can derive many interesting properties of the graph $G$ from the
value of $\lambda (G)$. The general phenomenon here is that if
$\lambda(G)$ is significantly less than $d$, then the edges of $G$
distribute like those of a  random graph with edge density $d/n$
\cite{AM, Thom, CGW}. One can use this  information  to derive
various properties of the graphs (see \cite{KSsur} for many results
of this kind). The whole concept can be generalized for non-regular
graphs. In this case, one needs to consider the Laplacian rather
than the adjacency matrix (see, for example, \cite{Chung}).

Estimating $\lambda(G)$ for a random regular graph is a well known
problem in the discrete math/theoretical computer science
community. A consequence of the well known Alon-Boppana  bound
\cite{AB} asserts that if $d$ is fixed and $n$ tends to infinity,
a.s.

$$\lambda_2(G_{n,d} )  \ge 2 \sqrt{d-1} - o(1). $$

\noindent Since $\lambda(G) \ge |\lambda_2 (G)|$  it follows that
a.s.

$$ \lambda(G_{n,d} )  \ge 2 \sqrt{d-1} - o(1). $$

Alon \cite{AB} conjectured that for any fixed $d$, a.s.

$$ \lambda_2 (G_{n,d}) = 2 \sqrt{d-1} +o(1). $$

\noindent Friedman \cite{Fri1} and Kahn and Szemer\'edi \cite{KSz}
showed that if $d$ is fixed and $n$ tends to infinity, then a.s.
$\lambda (G_{n,d}) = O(\sqrt d)$. Recently, Friedman, in a highly
technical paper \cite{Fri2}, proved Alon's conjecture. In fact, he
proved the stronger statement that a.s. $ \lambda (G_{n,d}) = 2
\sqrt{d-1} +o(1). $ This, together with the lower bound above,
determines the asymptotic of $\lambda(G_{n,d})$.

\begin{theorem} \label{theorem:Friedman} \cite{Fri2} (Second eigenvalue of random regular
graphs with fixed degree) For any fixed $d$ and $n$ tending to
infinity,  a.s.

$$\lambda (G_{n,d} ) = (2+o(1)) \sqrt {d-1}.$$

\end{theorem}

A $d$-regular graph $G$ is {\it Ramanujan} if $\lambda(G) \le 2
\sqrt {d-1}$. Explicit constructions of  Ramanujan graphs are highly
non-trivial and usually come from deep results in number theory (see
\cite{LPS} or \cite{Mar}, for example). On the other hand, the
following conjecture has been circulated in the last few years
(mentioned to the author by Sarnak)

\begin{conjecture} \label{conj:Sarnak} For $d$ fixed and $n$ tends
to infinity, $G_{n,d}$ is Ramanujan with positive constant
probability. \end{conjecture}

So far, we discussed the case when $d$ is a constant. What happens
if $d$ also tends to infinity with $n$? It is not clear (at least to
the author) that Friedman's  proof of Alon's conjecture in
\cite{Fri2} can be extended to this case.
 On the other hand, it is not hard to show that $\lambda (G(n,p))$, where $G(n,p)$ is
the Erd\H {o}s-R\'eyi random graph, is $(2+o(1)) \sqrt {np(1-p) }$
for sufficiently large $p$ (e.g., $p \ge n^{-1+\ep}$ for any fixed
$0 < \ep < 1$). Motivated by the universality principle, we make the
following conjecture

\begin{conjecture} \label{conjecture:normreg} Assume that $d \le n/2$  and both
$d$ and $n$ tend to infinity. Then a.s

$$\lambda (G_{n,d} ) = (2+o(1)) \sqrt {d(1-d/n)}. $$ \end{conjecture}

 Nilli \cite{Nil1} showed that for any  $d$-regular graph $G$ having two edges with
distance at least $2k+2$ between them $ \lambda_2 (G) \ge 2 \sqrt
{d-1} - 2 \sqrt {d-1}/(k+1)$. If $d= n^{o(1)}$ then $G_{n,d}$ has
diameter $\omega (1)$ with probability $1-o(1)$. Thus in this case

$$ \lambda (G_{n,d} ) \ge \lambda_2 (G_{n,d}) \ge  (2+o(1)) \sqrt {d} $$

\noindent with probability a.s. This proves the lower bound in
Conjecture \ref{conjecture:normreg}. For a general $d$,  it is easy
to show (by computing the trace of the square of the adjacency
matrix) that any $d$-regular graph $G$ on $n$ vertices satisfies

$$ \lambda (G) \ge \sqrt { d(n-d)/ (n-1)} \approx \sqrt {d(1- d/n)} . $$

(We would like to thank N. Alon for pointing out this bound.)

Let us now turn to the upper bound.  For $d= o(n^{1/2})$, one can
follow the Kahn-Szemer\'edi approach to show that $\lambda
(G_{n,d})= O(\sqrt d)$ a.s. For larger $d$, there is a weaker bound
$o(d)$ \cite[Theorem 2.8]{KSVW} proved by the trace method. The
following two approaches look promising:

\begin {itemize}

\item (Suggested by Krivelevich) Combine the sharp concentration
result in the previous section with the probability that a random
graph is regular. Using this, one can show for example that $\lambda
(G_{n,d}) = O( \sqrt {d \log n})$ for $d$ close to $n$ ($d= n/2$,
for instance).

\vskip2mm

\item The Sandwich Theorem  (Theorem \ref{theorem:sandwich})
implies $$\lambda (G_{n,d}) = \lambda(G(n,d/n)) + O( \sqrt {d \log
n}). $$

For most values of $d$, $\lambda(G(n, d/n)) = O(\sqrt {d})$. Thus,
if the Sandwich conjecture holds, it would imply a upper bound of
$O( \sqrt {d \log n})$ for most values of $d$.
\end{itemize}

The author feels confident that one can prove that $\lambda
(G_{n,d}) = O( \sqrt {d \log n})$ for all $d$ using these
approaches. However, removing the $\log$ term  seems non-trivial. In
fact, even the following special and weakened case  already looks
challenging

\vskip2mm

{\bf Problem.} {\it  Prove that  $\lambda (G_{n, n/2}) = O(\sqrt n)$
almost surely.}

\section {Determinant} \label{det}

The problem of determining the determinant of $M_n$ has been
considered by various researchers for at least 40 years. It was
proved by Koml\'os in 1967 \cite{Kom} that almost surely $\det M_n$
is not zero. In fact, it is easy to see that $\det M_n$ is divisible
by $2^{n-1}$, thus it follows that a.s. $|\det M_n| \ge 2^{n-1}$.
From above, Hadamard's inequality implies that $| \det M_n | \le
n^{n/2}$ (notice that all row vectors of $M_n$ have length $\sqrt
n$. It was often conjectured that with probability close to 1,
$|\det M_n|$ is close to this upper bound.

\begin{conjecture} \label{conj:determinant} Almost surely $| \det M_n | = n^{(1/2-o(1)) n}$.   \end{conjecture}

This conjecture is supported by the following observation of
Tur\'an, whose proof is a simple application of the linearity of
expectation.

\begin{fact} \label{fact:Turan}      $$\E ((\det M_n)^2) = n!. $$     \end{fact}

It follows immediately by Markov's bound that for any function
$\omega (n)$ tending to infinity with $n$, almost surely

$$|\det M_n | \le \omega (n)\sqrt {n!} . $$

\noindent Tao and Vu \cite{TVdet} established the matching lower
bound, which confirms Conjecture \ref{conj:determinant}.

\begin{theorem} \label{theorem:TVdet} Almost surely
$$|\det M_n |  \ge \sqrt {n!} \exp (- 29 \sqrt {n \log n}). $$ \end{theorem}

We are going to sketch the proof very briefly as it contains a
useful lemma. For a more detailed proof, we refer to \cite{TVdet}.

\begin{proof} We view $|\det M_n|$ as the volume of the
parallelepiped spanned by $n$ random $\{-1,1\}$ vectors. This volume
is the product of the distances from the $(d+1)$st vector to the
subspace spanned by the first $d$ vectors, where $d$ runs from $0$
to $n-1$. We are able to obtain  very tight control of this distance
(as a random variable), thanks to the following lemma.

\begin{lemma} \label{lemma:distance}
Let $W$ be a fixed subspace of dimension $1 \le d \le n-4$ and $X$
a random $\pm 1$ vector. Then
\begin{equation}\label{mean}
\E( \dist(X,W)^2 ) = n-d.
\end{equation}
Furthermore, for any $t >0$
\begin{equation}\label{pdw}
\P( |\dist(X, W) - \sqrt{n-d}| \ge t+1 ) \le  4 \exp(  -t^2/16).
\end{equation}
\end{lemma}

Observe that in this lemma, we do not need to assume that $W$ is
spanned by random vectors.  The lemma, however, is not applicable
when $d$ is very close to $n$ as it does not imply that the distance
is positive  almost surely. In this case, we do need to use the
assumption that $W$ is random. This assumption allows us to derive
information about the normal vector of $W$, which, combined with
Erd\H {o}s-Littlewood-Offord bound (see Theorem \ref{theorem:LO}),
provides control of the last few distances.
\end{proof}

\begin{remark}  After having written \cite{TVdet}, Tao and the author  discovered that Girko
(Section 6 of \cite{Girkobook}) claimed a very general theorem which
implies Theorem \ref{theorem:TVdet}. We are not able to understand
his proof and have not found anyone who does.
\end{remark}

Theorem \ref{theorem:TVdet} can be extended for much more general
models of random matrices.  Let $\xi_{ij}$, $1\le i,j \le n$, be
independent (but not necessarily i.i.d.) r.v's with the following
two properties:

\begin{itemize}

\item Each $\xi_{ij}$ has mean zero and variance one.

\item There is a constant $K$ that $|\xi_{ij}| \le K$ with
probability one.

\end{itemize}

\begin{theorem} \label{theorem:TVdetgeneral}
Consider the random matrix $M'_n $ with entries $\xi_{ij}$ as
above. Let $\ep$ be an arbitrary  positive constant. With
probability $1-o(1)$,

$$ |\det M'_n| \geq  \sqrt {n!} \exp( -n^{1/2+\ep}  ). $$

\end{theorem}

In certain situations, the assumption that $|\xi_{ij}|$ are
bounded from above by a constant is too strong. We are going to
consider the following less restricted model. Let $\xi_{ij}, 1\le
i,j \le n$ be i.i.d. random variables with mean zero and variance
one. Assume furthermore that their fourth moment is finite.
Consider the random matrix $M^{''}_n$ with $\xi_{ij}$ as its
entries.

\begin{theorem} \label{theorem:bai} \cite{TVdet}  We have,  with probability $1-o(1)$,
that

$$|\det M^{''}_n | \ge n^{(1/2-o(1)) n}. $$

\end{theorem}

An open problem concerning determinants  is to extend Theorem
\ref{theorem:TVdet} to symmetric matrices. Let $Q_n$ denote the
random symmetric matrix whose upper diagonal entries are i.i.d.
Bernoulli random variables.

\begin{conjecture} \label{conj:detsym} Almost surely, $|\det Q_n| = n^{(1/2+o(1))n }$. \end{conjecture}

It was proved only very recently \cite{CTV} that a.s $|\det Q_n|$ is
positive (which can be seen as the symmetric version of Koml\'os'
theorem mentioned above). The main difficulty here is that the row
vectors of $Q_n$ are no longer independent and so Lemma
\ref{lemma:distance} is not applicable as  there is a correlation
between the subspace $W$ and the vector $X$.

Finally, let us briefly discuss the situation with the permanent.
Notice that the estimate in Fact \ref{fact:Turan} is still valid for
the permanent of $M_n$. Thus, one would expect that, like the
determinant, the permanent of $M_n$ is typically of order
$n^{(1/2+o(1))n }$ (in absolute value). However, the following
problem is still open.

{\bf Problem.} {\it Show that the permanent of $M_n$ is a.s. not
zero. }

\section{Rank and singular probability: non-symmetric models}
\label{rank}

Let us consider the basic model $M_n$ and let  $p_n$ be the
probability that $M_n$ is singular. Estimating $p_n$ is well known
problem in discrete probability. From below it is clear that $p_n
\ge (1/2+o(1))^n$, as a matrix is singular if it has two equal
rows. A famous conjecture in the field asserts that this trivial
lower bound is sharp.

\begin{conjecture} \label{conjecture:singularity} $p_n =(1/2+o(1))^n $.\end{conjecture}

There is a  refined version of  the above conjecture where the right
hand side is more precise (see \cite{KKS}). However, Conjecture
\ref{conjecture:singularity}, as formulated,  is still open.

It is already non-trivial to show that $p_n=o(1)$. As mentioned in
the previous section, this was  done by Koml\'os almost fourty years
ago \cite{Kom}. The bound on $p_n$ in his original proof tends very
slowly to zero with $n$. Later, he found a new proof which showed
$p_n=O(n^{-1/2})$. In 1995, a breakthrough by Kahn, Koml\'os and
Szemer\'edi \cite{KKS} yielded  the first exponential bound $p_n =
O(.999^n)$. Their arguments were simplified by Tao and Vu in 2004
\cite{TVdet}, resulting in a slightly better bound $O(.958^n)$ and a
somewhat simpler proof. Shortly afterwards, Tao and Vu \cite{TVsing}
combined the approach from \cite{KKS} with ideas from additive
combinatorics to obtained the following  more significant
improvement

\begin{theorem} \label{theorem:TVsing} \cite{TVsing}  $p(n) \le (3/4+o(1))^n $.
\end{theorem}

The proof in \cite{TVsing} is highly technical and requires many
tools from discrete Fourier analysis, additive combinatorics and the
geometry of numbers. On the other hand, the proving scheme is
flexible and can be adapted to other models, sometime yielding
(surprisingly) sharp bounds. Let us present one such result. Instead
of $M_n$ we consider the following  ''lazy" model $M_n^{lazy}$. The
entries of $M_n^{lazy}$  are i.i.d random variables which equal zero
with probability one half and $1$ and $-1$ with probability one
quarter. (If one thinks  of the entries of $M_n$ as fair coin flips,
then in the ''lazy" model about half of the time we are lazy and
simply write zero instead flipping a coin.) It is clear that for the
lazy model the singular probability $p_n^{lazy}$ is again at least
$(1/2+o(1))^n$ (which is the probability that there is a zero row).
We are able to show that this bound is actually sharp

\begin{theorem} \label{theorem:TVrankmubounded} \cite{TVW} $p_n^{lazy}=(1/2+o(1))^n .$  \end{theorem}

Let us conclude this section by  two conjectures motivated by
studies from random graphs. These questions concern the resilience
of a structure, introduced in \cite{KimV, SudV}. Roughly speaking,
the resilience of a structure $S$ with respect to a property $\CP$
measures how much we have to change $S$ in order to destroy $\CP$.
We would like to measure the resilience of $M_n$ with respect to the
property of being non-sigular.

Given $\{-1,1\} $ matrix $M$, we denote by $Res (M)$ the minimum
number of entries we need to switch (from $1$ to $-1$ and vice
versa) in order to make $M$ singular. If $M$ is a sample of $M_n$,
it is easy to show that $Res (M)$ is, a.s, at most $(1/2+o(1))n$,
as we can, a.s, change that many entries in the first row to make
the first two rows equal. We conjecture that this is the best one
can do.

\begin{conjecture} \label{conjecture:effect1} Almost surely $Res (M_n)= (1/2+o(1))n $. \end{conjecture}

\noindent A closely related question (motivated by the notion of
local resilience from \cite{SudV}) is the following. Call a
$\{-1,1\}$ ($n$ by $n$) matrix $M$ {\it good} if all matrices
obtained by switching (from $1$ to $-1$ and vice versa) the
diagonal entries of $M$ are non-singular (there are $2^n$ such
matrices).

\begin{conjecture} \label{conjecture:effect2}  Almost surely $M_n$ is good.  \end{conjecture}

\section{Rank and singular probability: symmetric models}
\label{rank2}

Let us now consider symmetric matrices. The symmetric counterpart of
$M_n$ is $Q_n$. In fact is is more convenient to  consider
$Q(n,1/2)$ instead of $Q_n$, as the graph terminology is more
convenient and leads to natural extensions. (It is easy to show that
if  $Q(n,1/2)$ is a.s. non-singular then $Q_n$ is and vice versa.)

While the non-singularity of $M_n$ has been known for forty years
since \cite{Kom}, that of $G(n,1/2)$ was established only recently
by Costello, Tao and Vu \cite{CTV}. This confirmed a conjecture of
B. Weiss (this conjecture was communicated to the author by G.
Kalai and N. Linial).

\begin{theorem} \label{theorem:CTVrank} $Q(n,1/2)$ is a.s non-singular.  \end{theorem}

As  pointed out earlier in Section \ref{det}, the main difficulty in
going from the non-symmetric setting to the symmetric one is that
the row vectors in a random symmetric matrix are not independent.
The key tool that helped us to overcome this difficulty was the
so-called quadratic Littlewood-Offord inequality (Theorem
\ref{theorem:quadLO}), discussed in Section \ref{LO}.

It is natural to ask if Theorem \ref{theorem:CTVrank}  still holds
for a smaller density $p$. The answer is negative after a certain
threshold. Indeed, if $p < (1-\ep)\log n/n$ for some positive
constant $\ep$, then $G(n,p)$ has a.s. isolated vertices which
means that  its adjacency matrix has all zero rows and so is
singular. Costello and Vu proved that $\log n/n$ is the right
threshold.

\begin{theorem} \label{theorem:CVrank2} \cite{CV}  For any constant $\ep>0$,
$Q(n,(1+\ep)\log n/n)$ is a.s. non-singular.  \end{theorem}

 It remains  an interesting problem to estimate the rank of $Q(n, p)$ for $p < \log
 n/n$. The answer here is not yet conclusive, but some partial
 results are known. For instance, it was shown in \cite{CV} that if
 $p > (1+\ep) \log n/2n$, then the rank of $Q(n,p)$ is a.s. equal
  $n$ minus the number of isolated vertices in the graph. (The
 upper bound is trivial.)

Let us conclude  this section by stating two conjectures. The first
is a variant of Conjecture \ref{conjecture:singularity}. We denote
by $p^{sym}_n$ the probability that $Q(n,1/2)$ is singular. It is
easy to show that this probability is at least $(1/2+o(1))^n$.

\begin{conjecture} \label{conjecture:ranksym}  $p^{sym} _n = (1/2+o(1))^n .$  \end{conjecture}

This conjecture is perhaps very hard. The current best upper bound
on $p^{sym}_n$ is $n^{-1/4+o(1)}$. It seems already non-trivial to
replace $1/4$ by an arbitrary constant $C$.

The second conjecture concerns  random regular graphs. If $d=2$,
$G_{n,d}$ is a union of disjoint cycles and it is easy to show that
its adjacency matrix $Q_{n,d}$ is a.s. singular, as many of these
cycles have length divisible by 4. We conjecture that this is the
only case.

\begin{conjecture} \label{conjecture:rankreg} For all $d \ge 3$, $Q_{n,d}$ is a.s. non-singular. \end{conjecture}

\section{The condition number} \label{cond}

 For a matrix  $M$ the \emph{condition number} $c(M)$ is defined as
$$c(M):= \| (M) \|  \cdot  \| (M^{-1}) \|. $$
We adopt the convention that $c(M)$ is infinite if $M$ is not
invertible.

The condition number  plays a
 crucial role in applied linear algebra. In
particular, the complexity of any algorithm which requires solving
a system of linear equations usually involves the condition number
of a matrix \cite{BT}. Another area of mathematics where this
parameter is important is the theory of probability in Banach
spaces (see \cite{Rud} and the references therein).

The condition number of a random matrix is a well-studied object
(see \cite{Ede} and the references therein). In the case when the
entries of $M$ are i.i.d Gaussian random variables (with mean zero
and variance one), Edelman \cite{Ede} proved

\begin{theorem} \label{theo:condition1} Let $N_n$ be a $n \times n$ random matrix, whose
entries are i.i.d Gaussian random variables (with mean zero and
variance one). Then $\E(\ln c(N_n)) = \ln n + c+ o(1)$, where $c>0$
is an explicit constant.    \end{theorem}

In applications, it is usually useful to have  a tail estimate. It
was shown by Edelman and Sutton  \cite{ET} that

\begin{theorem} \label{theo:condition11} Let $N_n$ be a $n$ by $n$ random matrix, whose
entries are i.i.d Gaussian random variables (with mean zero and
variance one). Then for any constant $A >0$,
$$\P( c(N_n) \ge n^{A+1}) = O_A( n^{-A}) . $$
\end{theorem}

With the universality principle, it is reasonable to conjecture that
this estimate holds for the random Bernoulli matrix $M_n$ as well
(see \cite{ST} for an even more precise conjecture). However, this
seems very hard to prove. On the other hand, the following was
obtained recently by Tao and Vu \cite{TVinverse}

\begin{theorem} \label{theo:condition2} For any positive constant
$A$, there is a positive constant $B$ such that

$$\P( c(M_n) \ge n^{B}) \le n^{-A}. $$
\end{theorem}

It is well known that there is a constant $C$ such that the norm of
$M_n$ is at most $Cn^{1/2}$ with exponential probability $1- \exp(-
\Omega_\mu(n) )$ (in fact, one can prove this using the results in
Section \ref{norm}). Thus, Theorem \ref{theo:condition2} reduces to
the following lower tail estimate for the norm of  $M_n^{-1}$:

\begin{theorem} \label{theo:singular1}
For any positive constant $A$, there is a positive constant $B$ such
that

$$\P( \|M^{-1}_n \| \ge n^{B}) \le n^{-A}.$$

\end{theorem}

Shortly prior to Theorem \ref{theo:condition2}, Rudelson \cite{Rud}
proved the following result.

\begin{theorem} \label{theo:singular2}
 There are positive constants $c_1, c_2$ such that the following holds.
 For any $\ep \ge c_1 n^{-1/2}$
$$\P( \| M_n^{-1} \| \ge c_2 \ep n^{3/2}) \le \ep. $$
\end{theorem}

Both theorems can be generalized considerably (see \cite{TVinverse,
TVcond} and \cite{Rud}). Theorem \ref{theo:condition2} in particular
still holds if we replace $M_n$ by $M +M_n$ where $M$ is an
arbitrary matrix with polynomially bounded norm.

\begin{theorem} \label{theo:condition3} \cite{TVcond} For any positive constants
$A$ and $C$, there is a positive constant $B$ such that the
following holds. For any $n$ by $n$ matrix $M$ where $\|M \| \le
n^C$,

$$\P( c(M+ M_n) \ge n^{B}) \le n^{-A}. $$
\end{theorem}

The point here is that $M$ itself can have very large condition
number. (In fact if $M$ is singular then its condition number is
infinity.) Theorem \ref{theo:condition3} asserts that a Bernoulli
perturbation of $M$ has small condition number with high
probability. The Gaussian version of Theorem \ref{theo:condition3}
was proved by Spielman and Teng in \cite{ST}. For the connection
of these theorems to numerical  linear algebra and theoretical
computer science, we refer to \cite{ST} and \cite{TVcond}.

\section{Littlewood-Offord and quadratic Littlewood-Offord} \label{LO}

Let $\bv =\{v_1, \dots, v_n\}$ be a set of $n$ integers and $\xi_1,
\dots, \xi_n$ be i.i.d random Bernoulli variables. Define $S:=
\sum_{i=1}^n \xi_i v_i$ and $p_{\bv} (a) := \P (S=a)$ and $p_{\bv}
:= \sup_{a \in \BZ} p_{\bv } (a)$.

  Erd\H {o}s, answering a question of
 Littlewood and Offord, proved the following theorem, which we are referring
  to as the Erd\H {o}s-Littlewood-Offord inequality.

\begin{theorem} \label{theorem:LO} Let $v_1, \dots, v_n$ be non-zero
numbers and $\xi_i$ be i.i.d Bernoulli random variables. Then

$$ p_{\bv} \le \frac{{n \choose
{\lfloor n/2 \rfloor}}}{2^n} = O(\sqrt n). $$\end{theorem}

Theorem \ref{theorem:LO} is a classical result in combinatorics and
has many non-trivial extensions (see \cite[Chapter 7]{TVbook} or
\cite{Hal} and the references therein).

The random sum $S:= \sum_{i=1}^n \xi_i v_i$ plays a central role in
the study of random Bernoulli matrices. In many problems (such as
that of the determinant, rank or condition number), the critical
parameter is the distance from a random Bernoulli vector to the
hyperplane spanned by another $n-1$ random Bernoulli vectors (these
are the row vectors of $M_n$). If $(v_1, \dots , v_n)$ is the (unit)
normal vector of the hyperplane, then this distance is exactly the
absolute value of the random sum $S= \sum_{i=1}^n \xi_i v_i$.

Bounding the above mentioned distance  is  one of the main
difficulties when one goes from Gaussian matrices to Bernoulli
matrices. If the entries of the matrix are Gaussian, then the
distance in question is a simple object. Thanks to symmetry, the
position of the hyperplane does not really matter and so we can
condition on it. Furthermore, the distribution of the distance from
a random Guassian vector to a fixed hyperplane is well understood.
The situation in the Bernoulli case is very different. In this case,
the random vectors are chosen from the vertices of the
$n$-dimensional $\{-1,1\}$ cube. Very little is known about the
hyperplanes spanned by $n-1$ such vectors. Let us point out,
however, that there are planes that contain a constant fraction of
the vertices of the $\{-1,1\}$ cube and in this case the distance in
question is zero with constant probability.

In order to treat random symmetric matrices (in particular $Q(n,
1/2)$), Costello, Tao and Vu \cite{CTV} introduced the following
quadratic version of Theorem \ref{theorem:LO}

\begin{theorem} \label{theorem:quadLO} Let $c_{ij}, 1\le i, j \le n$ be non-zero numbers.
Consider the quadratic form $F= \sum_{1 \le i \le j \le n} c_{ij}
\xi_i \xi_j$, where $\xi_i, 1\le i\le n$ are i.i.d Bernoulli random
variables. Then for any $a$

$$\P (F=a) = O(n^{-1/8}). $$
\end{theorem}

The exponent $-1/8$ was improved to $-1/4$ (see \cite{CV}). It is
conjectured \cite{CV} that the sharp exponent would be $-1/2$. The
lower bound is given by the quadratic form $F = (\sum_{i=1}^n
\xi_i)^2$.

\section{Random walks and Lazy Random Walks} \label{randomwalk}

Let $\bv = \{v_1, v_1, \dots, v_n \}$ be the set of $n$ non-zero
numbers and consider the random walk $W$ on the real line (starting
at $0$) which at step $i$ goes to the left by $v_i$ with probability
one half and to the right by $v_i$ with probability one half.   The
probability $p_{\bv} (0)= \P(\sum_{i=1}^n \xi_i v_i=0)$ is exactly
the probability that $W$ returns to the origin after $n$ steps.

Let $\mu$ be a constant between $0$ and $1$ and consider the ''lazy"
random walk $W^{\mu}$ (starting at $0$) which at step $i$ stays with
probability $1-\mu$ and goes to the left by $v_i$ with probability
$\mu/2$ and to the right by $v_i$ with probability $\mu/2$. Let
$p_{\bv}^{\mu} (0)$ be the probability that the lazy walk return to
zero after $n$ steps.

Intuitively, one would expect that $p_{\bv}^{\mu}(0)$ is larger than
$p_{\bv}(0)$ (especially when $\mu$ is small), as the lazy walk has
a stronger tendency to stay near the starting point. Quantitatively,
one can show (using Fourier analysis and the elementary fact that
$|\cos x| \le 3/4 + 1/4 \cos 2x$) that for any $\bv$

$$p_{\bv}(0) \le p_{\bv} ^{1/4} (0). $$

The next question is: Can one improve this to

\begin{equation} \label{equ:lazy} p_{\bv} (0) \le  \ep p_{\bv} ^{1/4} (0),
\end{equation}

\noindent for any positive constant $\ep$ ? The answer is negative.
If we take $\bv= \{1,1, \dots, 1\}$, then it is easy to show that

$$p_{\bv} (0) = (c+o(1)) p_{\bv} ^{1/4} (0)  $$

\noindent where $c$ is a  positive constant (depends on $1/4$).
However, it is possible to classify all sets $\bv$ where
\eqref{equ:lazy} fails (under some slight assumptions). This
classification is the heart of the proof of Theorem
\ref{theorem:TVsing}. The precise statement is somewhat technical
(we refer to \cite{TVsing} for details), but roughly it says  that
if \eqref{equ:lazy} fails then $\bv$ is contained in a generalized
arithmetic progression  with constant rank and small volume.

\section{Inverse Littlewood-Offord theorems} \label{inverseLO}

 A
set $$P= \{c+ m_1a_1 + \dots +m_d a_d| M_i \le m_i \le M_i'\}$$
 is called a {\it generalized arithmetic
progression } (GAP) of rank $d$. It is convenient to think of $P$ as
the image of an integer box $B:= \{(m_1, \dots, m_d)| M_i \le m_i
\le M_i' \} $ in $\Z^d$ under the linear map
$$\Phi: (m_1,\dots, m_d) \mapsto c+ m_1a_1 + \dots + m_d a_d. $$ The
numbers $a_i$ are the {\it generators } of $P$.  For a set $A$ of
reals and a positive integer $k$, we define the iterated sumset
$$kA := \{a_1+\dots +a_k |a_i \in A  \}.$$

Let us take another look at Theorem \ref{theorem:LO}. This theorem
is sharp, as is shown by taking $v_1=v_2=\dots= v_n =1$. However,
the bound changes significantly if one forbid this special case.
Erd\H {o}s and Moser \cite{EM} showed that under the stronger
assumption that the $v_i$ are non-zero and different,

$$p_{\bv}   =O(n^{-3/2} \ln n). $$

 They conjectured that the logarithmic term is
not necessary and this was confirmed by S\'ark\"ozy and
Szemer\'edi \cite{SS} (see also \cite{Hal}). Again, the bound is
sharp (up to a constant factor), as can be seen by taking
$v_1,\ldots,v_n$ to be a proper arithmetic progression such as
$1,\ldots,n$.

In the above two examples, we observe  that in order to make
$p_{\bv}$ large, we have to impose a very strong additive
structure on $\Bv$ (in one case we set the $v_i$'s to be the same,
while in the other we set them to be elements of an arithmetic
progression). A more general example is the following

{\it Example.} Let $P$ be a GAP of rank $d$ and volume $V$. Let
$v_1, \dots, v_n$ be (not necessarily different) elements of $P$.
Then the random variable $S =\sum_{i=1}^n \xi_i v_i$ takes values in
the GAP $nP$ which has volume $n^d V$. From the pigeonhole principle
and the definition of $p_{\bv}$, it follows that
$$p_{\bv}  \ge n^{-d} V^{-1}.$$

\noindent This example shows that if the elements of $\Bv$ belong
to a GAP with small rank and small volume then  $p_{\bv}$ is
large. One may conjecture that  the inverse  also holds, namely,

\vskip2mm

{\it If $p_{\bv}$ is large, then (most of) the elements of $\Bv$
belong to a GAP with small rank and small volume. }

 Tao and Vu \cite{TVinverse} have managed to quantify this
statement.

\begin{theorem}\label{theorem:inverse3} \cite{TVinverse}
Let  $A, \alpha > 0$ be arbitrary constants. There are constants
$d$ and $B$ depending on $A$ and $\alpha$ such that the following
holds.
 Assume that $\Bv = \{v_1, \ldots, v_n\}$ is a multiset of integers satisfying
$\P (S=0) \geq n^{-A}$. Then there is a GAP $Q$ of rank at most
$d$ and volume at most $n^B$ which contains all but at most
$n^{\alpha}$ elements of $\Bv$ (counting multiplicity).

\end{theorem}

Notice that the small set of exceptional elements is not
avoidable. For instance, one can add $O(\log n)$ completely
arbitrary elements to $\Bv$,
 and only decrease $p_{\bv}$ by a factor of $n^{-O(1)}$ at worst.

Theorem \ref{theorem:inverse3} is one of the main ingredients in
the proofs of Theorems \ref{theo:condition2} and
\ref{theo:condition3}. For many other theorems of this type, see
\cite{TVinverse}.

{\bf Acknowledgement.} The author would like to thank K. Costello
and P. Wood for proofreading the manuscript.

\end{document}